\documentclass[11pt]{article}
\usepackage{latexsym}
\usepackage{theorem}
\usepackage{graphicx}
\usepackage{amsmath}
\usepackage{color}
\usepackage{xcolor}
\usepackage{amsfonts}
\usepackage{natbib}
\usepackage{soul}
\usepackage{hyperref}
\usepackage{enumitem}

\theorembodyfont{\rmfamily}
\newtheorem{theorem}{Theorem}

\newtheorem{lemma}[theorem]{Lemma}

\newtheorem{corollary}[theorem]{Corollary}

\theoremstyle{break}

\newcommand{\la}{\lambda}

\newcommand{\red}[1]{{\color{red}#1}}

\title{A Booby Trap Game}

\author{Thomas Lidbetter\thanks{Department of Management Science and Information Systems, Rutgers Business School, Newark, NJ 07102, tlidbetter@business.rutgers.edu} \and Kyle Y. Lin\thanks{Operations Research Department, Naval Postgraduate School, Monterey, CA 93943, kylin@nps.edu}}

\date{November 27, 2024}

\begin{document}

\maketitle

\begin{abstract}
\noindent
This paper presents a booby trap game played between a defender and an attacker on a search space, which may be a compact subset of Euclidean space or a network.
The defender has several booby traps and chooses where to plant them.
The attacker, aware of the presence of these booby traps but not their locations, chooses a subset of the space and collects a reward equal to the measure of the subset. 
If the attacker does not encounter any booby traps, then the attacker keeps the reward; otherwise, the attacker gets nothing.
The attacker's objective is to maximize the expected reward, while the defender's objective is to minimize it.
We solve this game in the case that the search space is a compact subset of Euclidean space, and then turn our attention to the case where the search space is a network in which the attacker must choose a \textit{connected} subset of the network.
We solve the game when the network is a circle or a line.
For the case of one booby trap, we solve the game for 2-connected networks, and when the network is a tree we present an upper bound and a lower bound for the value of the game whose ratio is at most $27/25$.
We also present an optimal solution for each player in a few cases where the tree is a star network.
\end{abstract}



\section{Introduction}

We consider an attacker-defender game in which a defender places a fixed number of booby traps in some search space and an attacker attempts to capture some subset of the space. If the chosen subset contains a booby trap, the attacker receives nothing; otherwise he receives a payoff equal to the length, area, or measure of the subset.
For example, the defender could represent a security team using hidden cameras to provide surveillance in an area, and the booby traps are security cameras.
The space could represent a house, a shop, a museum, an airport, or a corporate office.
The attacker could be a thief trying to steal valuable items, a drug smuggler trying to put illicit drug into luggage of random passengers, or a spy trying to plant their own eavesdropping devices, without getting caught on the security cameras.

\cite{lidbetter2020search} presents a game that shares some characteristics with the games analyzed in this paper. 
There are several boxes, each of which contains a reward.
One player hides a number of booby traps among some of these boxes, and the other player chooses which boxes to open to collect as much reward as possible without opening any booby-trapped boxes.
While such a formulation is reasonable in some cases where the search space can be represented by discrete locations, it is not applicable to search spaces such as fields, roads, paths, airports, where the topography is important.
In this paper, we study booby trap games played in a continuous search space.

\begin{figure}[ht]
\center
\includegraphics[width=10cm]{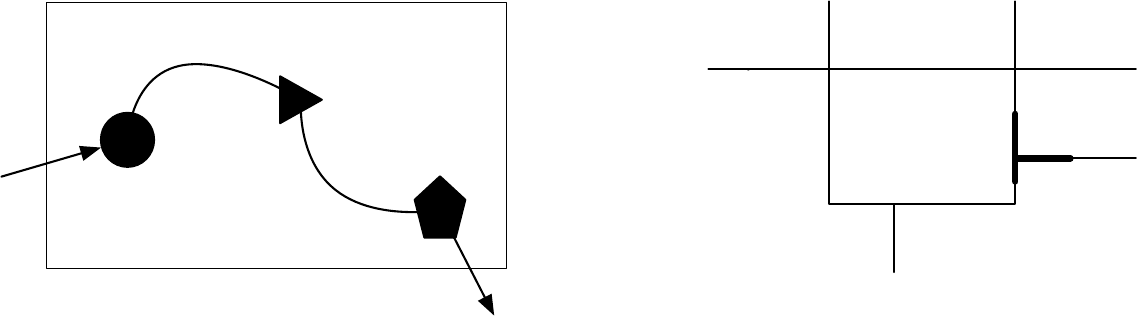}
\caption{The booby trap game played on a compact subset of a Euclidan space (left), and on a network (right).}
\label{fig:space}
\end{figure}

We consider two mathematical formulations of the search space in this paper. 
The first formulation deals with search spaces that are compact subsets of a Euclidean space, where the attacker can choose any subset of the space---such as the dark areas in the rectangle shown on the left-hand side of Figure~\ref{fig:space}.
For example, a farmer could place booby traps on a farm growing vegetables or fruits in order to catch rats and mice.
A retreating troop could place land mines to disrupt the operations of the opposing force.
The second formulation deals with networks, where the attacker chooses a connected subset of the network---such as the thickened subnetwork shown on the right-hand side of  Figure~\ref{fig:space}.
For example, a security team may install hidden surveillance cameras along hallways in a museum or corridors in a corporate firm trying to catch thieves or spies.

While the mathematical study of booby trap games is relatively new, our work is related to a broader research area known as search games.
In a typical search game, one player hides \textit{valuable} objects and the other player wants to \textit{find} those objects, whereas in booby trap games, one player hides \textit{harmful} objects and the other player wants to \textit{avoid} those objects.
Recent works on search games played in a discrete space include \cite{bui2024optimal}, \cite{clarkson2023classical}, \cite{clarkson2024computing}, \cite{duvocelle2022competitive}, \cite{lidbetter}, \cite{lidbetter2019searching}, and \cite{yolmeh2021weighted}.
Recent works on search games played in a continuous space include \cite{alpern2019stochastic}, \cite{alpern2013mining} \cite{angelopoulos2020competitive}, \cite{garrec2019continuous}, and \cite{dagan2008network}.
For an overview and review of search games, see \cite{AlpernGal}, \cite{Gal2011}, and \cite{Hohzaki}.

The rest of this paper proceeds as follows.
Section~\ref{sec:general} presents a complete solution of the booby trap game played on a compact subset of Euclidean space.
Section~\ref{sec:network} defines the booby trap game played on a network and solves the game for two important classes of networks: cycles and paths.
The remainder of the paper focuses on the case where the defender has one booby trap, and we conclude Section~\ref{sec:network} with a solution for 2-connected networks.
Section~\ref{sec:tree} concerns tree networks and presents strategies for both players that guarantee an expected payoff that is within $27/25$ of the value of the game.
Section~\ref{sec:star} focuses on star networks, and solves the game for symmetric stars and stars with three arcs.
Section~\ref{sec:conclusions} concludes.

\section{Booby Trap Game on a Subset of Euclidean Space} 
\label{sec:general}
In this section we consider a game played on a compact subset $Q$ of Euclidean space, equipped with Lebesgue measure $\lambda(\cdot)$.
For example, $Q$ might be a 2-dimensional shape representing a field, so that Lebesgue measure corresponds to area. By rescaling, we may assume that the measure of $Q$ is $1$. A defender has $k$ booby traps and chooses where to plant them in $Q$, so that a pure strategy for the defender is any $k$ distinct points in $Q$. 
The attacker, aware of the presence of the $k$ booby traps but not their locations, chooses a measurable subset $S$ of $Q$.
If no booby trap is planted in $S$, then the attacker earns reward $\lambda(S)$; otherwise, the attacker gets nothing.
The attacker's objective is to maximize the expected reward, while the defender's objective is to minimize it. We denote this zero-sum game by $\Gamma_k(Q)$.

For example, the space could represent a farm growing carrots or potatoes.
The attacker is a mouse and the defender is a farmer planting mousetraps.
For another example, the space could represent a battlefield left behind by the retreating troop.
The defender is the retreating troop who planted landmines in the field, and the attacker is their enemy who may want to use the field to set up their temporary military base.


\begin{theorem}
For a compact subset $Q$ of Euclidean space, the game $\Gamma_k(Q)$ has a value $k^k / (k+1)^{k+1}$. It is optimal for the defender to place each booby trap uniformly randomly on $Q$. The attacker has strategies that guarantee a payoff that is arbitrarily close to the value of the game.
\end{theorem}
\textit{Proof.}
Consider the defender strategy that places each booby trap uniformly randomly in $Q$, independent of one another.
If the attacker chooses a subset with measure $x$, then the expected reward is $x (1-x)^k$.
Calculus shows that $x = 1/ (1+k)$ maximizes the preceding and the maximal value is $k^k / (k+1)^{k+1}$, which is therefore an upper bound for the value.

Next consider a strategy for the attacker.
Partition $Q$ into $m \ge k+1$ subsets, each having measure $1/m$.
If the attacker chooses $r \le m$ of these subsets uniformly randomly, then the probability that none of these $r$ subsets contains a booby trap is at least
\[
\frac{\binom{m-k}{r}}{\binom{m}{r}} = \frac{(m-k)! \, (m-r)!}{(m-k-r)! \, m!} = \prod_{i=0}^{k-1} \left( 1 - \frac{r}{m-i} \right),
\]
which is achieved if each booby trap is in a different subset.
The expected reward is therefore at least
\[
\frac{r}{m} \prod_{i=0}^{k-1} \left( 1 - \frac{r}{m-i} \right).
\]
Take $r = \lceil m / (k+1) \rceil$, so 
\[
\frac{m}{k+1} \le r < \frac{m}{k+1} + 1.
\]
We can get a lower bound for the expected reward with
\[
\frac{r}{m} \prod_{i=0}^{k-1} \left( 1 - \frac{r}{m-i} \right) \ge \frac{1}{k+1} \prod_{i=0}^{k-1} \left( 1 - \frac{m+k+1}{m-i} \frac{1}{k+1} \right).
\]
By taking $m \rightarrow \infty$, the right-hand side of the preceding converges to
\[
\frac{1}{k+1} \left( 1 - \frac{1}{k+1} \right)^k = \frac{k^k}{(k+1)^{k+1}}.
\]
It follows that the attacker can guarantee an expected reward arbitrarily close to $k^k / (k+1)^{k+1}$, which is therefore a lower bound for the value.
Consequently, we prove the value of the game as stated.
\hfill $\Box$

It is interesting to point out that by writing the value of the game as $(1/k)(1-1/(k+1))^{k+1}$, it is clear that as $k \rightarrow \infty$, the value is asymptotically equal to $1/(ke)$.

\section{Booby Trap Game on a Network}
\label{sec:network}

We now consider an alternative version of the game introduced in the previous section played on a network. To define precisely what we mean by a network, we start with a graph $G$ with edges and vertices, and assume that every edge $e$ is assigned a length $\lambda(e)$. The edge $e$ is then identified with an open interval of length $\lambda(e)$, endowed with Lebesgue measure and Euclidean distance $d$. Thus, we consider $\lambda$ as a measure (which we refer to as {\em length}) on the whole network $Q$. Also, $d$ naturally extends to a metric on $Q$, given by the shortest path length between two points. Vertices and edges in $G$ correspond to what we refer to as {\em nodes} and {\em arcs} in $Q$. We make a standing assumption that, by rescaling, the length of the network is always equal to $1$.

As before, the defender chooses $k$ points to place booby traps, but this time the attacker chooses a {\em connected} subset of $Q$.
If the subset contains no booby trap, then the attacker earns a reward equal to the length of the subset; otherwise, the attacker gets 0. We denote this game $\Gamma'_k(Q)$.
We first solve the game $\Gamma'_k$ played on a circle.

\begin{theorem}
Consider the game $\Gamma'_k(C)$ played on a circle $C$.
An optimal strategy for the attacker is to choose a subinterval of length $1/ (2k)$ uniformly randomly on the circle.
An optimal strategy for the defender is to place the $k$ booby traps at equal distances on $C$.
The value of the game is $1/ (4k)$.
\end{theorem}
\textit{Proof.}
Suppose that the attacker selects a subinterval of length $1/ (2k)$ uniformly randomly on the circle.
At least half of these intervals do not contain a booby trap, so the expected payoff is at least $(1/2)\cdot 1/(2k) = 1/(4k)$, which is a lower bound for the value of the game.

Suppose the defender places $k$ booby traps equal distance on the circle, so the adjacent booby traps are $1/k$ apart.
If the attacker chooses a subinterval of length $1/k$ or longer, then the subinterval contains at least 1 booby trap, so the reward is 0.
If the attacker chooses a subinterval of length $x < 1/k$, the probability that it contains a booby trap is $k x$, so the expected payoff is $x(1-k x)$, which is maximized at $x=1/(2k)$, giving an expected payoff of $1/(4k)$, which is an upper bound for the value of the game.

Because the attacker can guarantee expected reward at least $1/(4k)$ and the defender can guarantee expected reward at most $1/(4k)$, the value of the game is $1/(4k)$.
\hfill $\Box$

\bigskip

The idea used to solve the game played on a circle can be used to solve the game played on a line segment, as presented below.

\begin{corollary}
Consider the game $\Gamma'_k([0,1])$ played on a line segment $[0,1]$.
An optimal strategy for the attacker is to divide $[0,1]$ into $2k$ subinterval each having length $1/(2k)$, and choose one subinterval uniformly randomly.
An optimal strategy for the defender is to select a number $u$ uniformly randomly over $(0, 1/k)$, and place the $k$ booby traps at $u, u+1/k, u + 2/k, \ldots, u+ (k-1)/k$.
The value of the game is $1/ (4k)$.
\end{corollary}

For networks other than circles and lines, each player's strategy is strongly dependent on the topology of the network.
If the defender has $k \ge 2$ booby traps, then generally speaking, the defender wants to randomize their placements but also keep the booby traps apart to cover all corners the network evenly.
If the defender has $k=1$ booby trap, then it is possible to solve the game for networks with certain properties.
We begin with the following lemma.

\begin{lemma} \label{lem:k=1}
Consider the game $\Gamma_1(Q)$. The value of the game is at most $1/4$, with equality if there exists a partition of $Q$ into two connected subsets each having length $1/2$. In this case, the uniform strategy is optimal for the defender.
An optimal strategy for the attacker is to take each subset with probability $1/2$.
\end{lemma}
\textit{Proof.}
Suppose the defender places the booby trap uniformly randomly on the network.
If the attacker chooses a subnetwork of length $x$, then the expected reward is $x (1-x)$, which is maximized for $x = 1/2$, giving an expected reward of $1/4$, which is an upper bound for the value of the game.

If there is a partition of the network into two connected subsets each having length $1/2$, then the attacker can choose each subset with probability $1/2$ to ensure an expected reward $1/4$, so the value of the game is $1/4$ in this case.
\hfill $\Box$

\bigskip

A network is called {\em 2-connected} if it remains connected upon removal of any node.
\cite{lempel1967algorithm} showed that given any arc with endpoints $s$ and $t$ in a 2-connected network with $n$ nodes, the nodes of the network may be numbered from $1$ to $n$ so that node $s$ receives number $1$, node $t$ receives number $n$, and every other node is adjacent both to a lower-numbered and to a higher-numbered node. 
Such a labeling is called an {\em $st$-numbering}, and \cite{lempel1967algorithm} presented an algorithm for computing an $st$-numbering. 
A more efficient algorithm was given in \cite{even1976computing}.

Using the concept of an $st$-numbering, we show in the following lemma that any 2-connected network can be partitioned into two connected subnetworks of equal length.
For a subset $A$ of nodes of a network, we call the subnetwork that includes all the nodes of $A$ and all arcs whose endpoints are both in $A$ the \textit{subnetwork induced by $A$}.

\begin{lemma} \label{lem:partition}
Any 2-connected network with total length 1 can be partitioned into two subnetworks each having length $1/2$.
\end{lemma}
\textit{Proof.}
Let $s$ and $t$ be the endpoints of an arbitrary arc of $Q$, and suppose an $st$-numbering is given. For $j=1,\ldots,n-1$, let $Q_j$ be the subnetwork induced by the nodes with labels $1,2,\ldots,j$. It is easy to prove by induction on $j$ that $Q_j$ is connected. Indeed, it is trivially true for $j=1$; for $j \ge 2$, assuming the claim is true for all smaller $j$, the $st$-numbering ensures that the node with label $j$ is adjacent to some node of $Q_{j-1}$, which is connected, by the induction hypothesis.

Let $\overline{Q_j}$ denote the complement of $Q_j$. The subnetwork $\overline{Q_j}$ is also connected, by a similar argument applied to the subnetwork induced by the nodes with labels $j+1,\ldots,n$, along with the interiors of all arcs with one endpoint in $Q_j$ and one in $\overline{Q_j}$.

The sequence $\lambda(Q_1),\ldots,\lambda(Q_n)$ is increasing, with $\lambda(Q_1)=0$ and $\lambda(Q_n)=1$. If there is some $j$ with $\lambda(Q_j) = 1/2$, then the proof is complete, since $Q$ can be partitioned into $Q_j$ and $\overline{Q_j}$, which both have length $1/2$.

If $\lambda(Q_j) \neq 1/2$ for $j=1,\ldots,n$, then there must exist some $j$ such that $\lambda(Q_j) < 1/2$ and $\lambda(Q_{j+1}) > 1/2$. In this case, let $\mu=1/2-\lambda(Q_j)$. Note that $Q_{j+1} \setminus Q_j$ consists of all arcs one of whose endpoints is the node labeled $j+1$ and the other is in $Q_j$. The total length of all these arcs is $\lambda(Q_{j+1})-\lambda(Q_j) > \mu$. We define a subnetwork $R$ of $Q$ recursively, starting with $Q_j$. One by one, in an arbitrary order, we add arcs in $Q_{j+1} \setminus Q_j$ to $R$ until reaching some arc $a$, whose addition would cause the length of $R$ to exceed $1/2$; at this point, we add a subinterval of $a$ chosen with one endpoint in $Q_j$, such that the length of $R$ reaches precisely $1/2$.  
The subnetwork $R$ and its complement are both connected and each has length $1/2$.
\hfill $\Box$

\bigskip

\begin{figure}[!ht]
\center
\includegraphics[width=7cm]{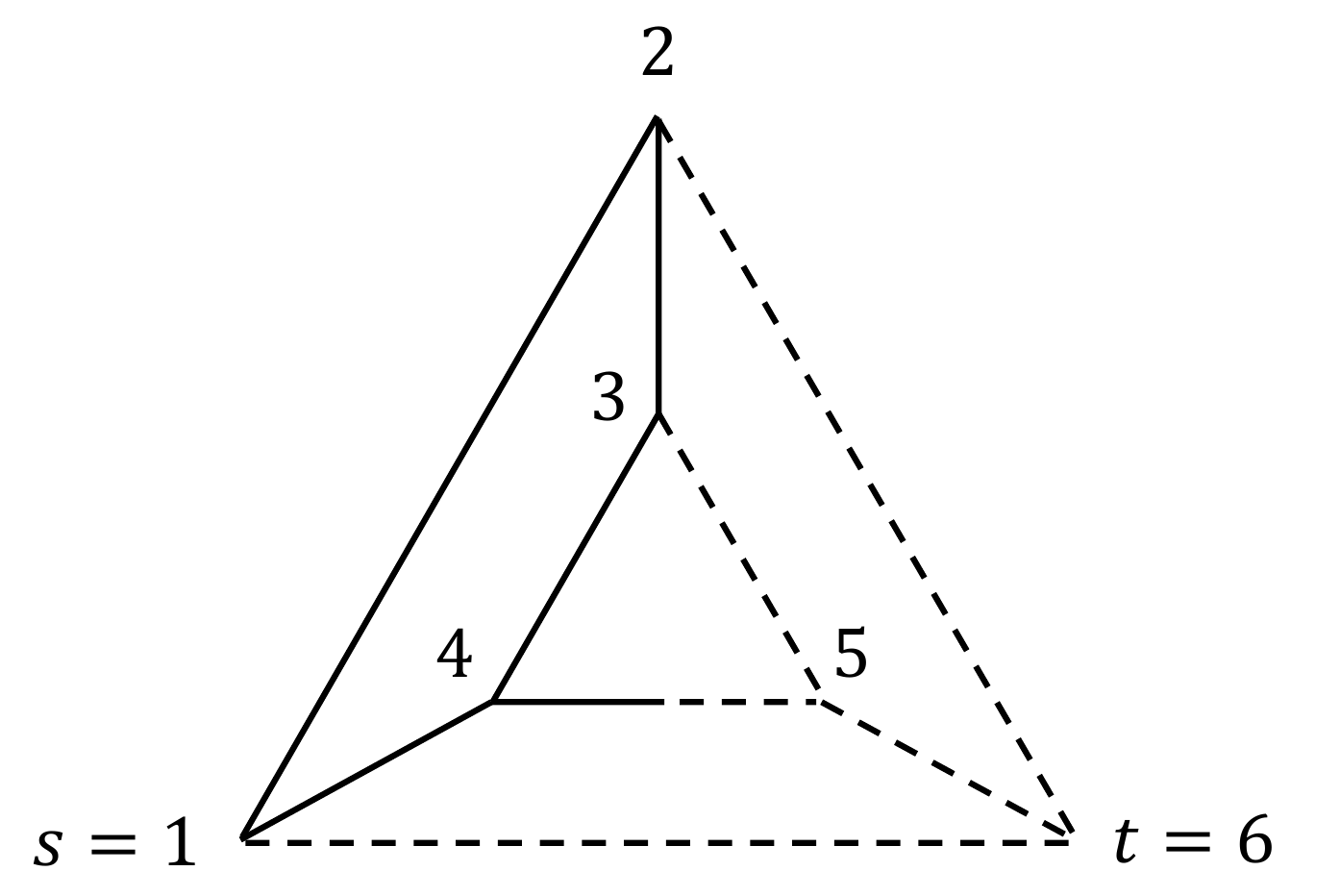}
\caption{The partition of a 2-connected network into two connected networks of equal length.}
\label{fig:partition}
\end{figure}

In Figure~\ref{fig:partition}, we illustrate the process of using an $st$-numbering to partition a 2-connected network into two subnetworks of equal length. For simplicity, suppose all the arcs in the network depicted have the same length of $1/9$. The nodes have been labeled with an $st$-numbering. The subnetwork $Q_4$ has length $4/9$ and $Q_5$ has length $6/9$. By adding half of the arc between nodes labeled $4$ and $5$ to $Q_4$, we obtain a connected subnetwork of length $1/2$ (shown with solid lines in the figure), whose complement (shown with dashed lines) is also connected and has length $1/2$.

\begin{theorem}
   If $Q$ is a 2-connected network, then the value of the game $\Gamma'_1(Q)$ is $1/4$.
\end{theorem}
\textit{Proof.}
    By Lemma~\ref{lem:partition}, the network $Q$ can be partitioned into two subnetworks of length $1/2$. By Lemma~\ref{lem:k=1}, the value of the game on $Q$ is $1/4$ and the defender's uniform strategy is optimal.
\hfill $\Box$

\bigskip

\section{Booby Trap Game on a Tree}
\label{sec:tree}
This section concerns the game $\Gamma'_1$ played on a tree with $k=1$ booby trap.
We first define the \textit{centroid} of a tree, and then use it to develop the \textit{centroid strategy} for the defender and the \textit{partition strategy} for the attacker.
The defender's centroid strategy provides an upper bound for the value of the game, and the attacker's partition strategy provides a lower bound.
These two bounds are always within 8\% of each other.
Finally, we extend the definition of the centroid to a general network, when appropriate.
While not all networks have a centroid, our results apply to those networks that have a centroid.

\subsection{Centroid of a Tree}
\label{sec:centroid}
For any point $x$ of a tree $Q$, let $Q(x)$ denote the network obtained by removing $x$. Unless $x$ is a leaf node, $Q(x)$ will be disconnected. 
Let $h(x)$ be the maximum length of any connected component of $Q(x)$. Observe that $h(x)$ is a lower semicontinuous function on a compact space (with respect to the shortest path metric). Therefore, it attains its minimum. We call any point $x^*$ that minimizes $h(x)$ a {\em centroid} of the tree and the minimized value $h(x^*)$ the \textit{radius} of the tree. 
This definition of the centroid is analagous to that of the centroid for a discrete graph; see \cite{slater1978centers}, for example.

\begin{lemma} \label{lem:centroid}
Any tree has a unique centroid, and its radius is at most $1/2$.
\end{lemma}
\textit{Proof.}
Let $x^*$ be a centroid of a tree $Q$, and suppose, for a contradiction, that $h(x^*) > 1/2$, which means that the largest connected component $A$ of $Q(x^*)$ has length greater than $1/2$. Let $x'$ be a point in $A$ at distance $\varepsilon < 2\lambda(A)-1$ from $x^*$ such that there are no nodes on the path from $x'$ to $x^*$, except possibly $x^*$. (In particular, $x'$ may not be a node.) 
The two connected components of $Q(x')$ have respective lengths $\lambda(A)-\varepsilon < \lambda(A)$ and $1-\lambda(A)+\varepsilon < \lambda(A)$, where the inequalities follow from the definition of $\varepsilon$. 
Consequently, $h(x') < h(x^*)$, so $x^*$ cannot be a centroid, which draws a contradiction.
We can therefore conclude that $h(x^*) \le 1/2$.

To show that the centroid $x^*$ is unique, we prove a stronger statement: There cannot be any other point $y \in Q$ with $h(y) \le 1/2$. Suppose such a point existed, and that the distance from $x^*$ to $y$ is $d >0$. Let $A_1$ be the component of $Q(x^*)$ containing $y$ and let $A_2$ be the component of $Q(y)$ containing $x^*$. It is easy to see that $A_1 \cup A_2 =Q$, so we must have
\[
1=\lambda(Q) = \lambda(A_1)+\lambda(A_2)- \la(A_1 \cap A_2) \le \frac{1}{2} + \frac{1}{2} -d <1,
\]
a contradiction.
Therefore, a tree has a unique centroid.
\hfill $\Box$

\bigskip

Locating the centroid of a tree is straightforward. First, if we find some node $x$ with $h(x) \le 1/2$, then $x$ must be the centroid and $h(x)$ the radius.
If none of the nodes is a centroid, then the centroid must be on the interior of the two adjacent nodes $y$ and $z$ such that $y$ lies in the largest component of $Q(z)$ and $z$ lies in the largest component of $Q(y)$.
In this case, the radius of the tree must be $1/2$.

If the radius of a tree is $1/2$, then the value of the game is $1/4$ according to Lemma~\ref{lem:k=1}.
The rest of this section concerns the case where the radius is less than $1/2$.
Because the radius is less than $1/2$, the centroid must disconnect the tree into at least 3 components.
Write $n$ for the number of components, and $a_1, a_2, \ldots, a_n$ for the lengths of these $n$ components.
Without loss of generality, label the components such that $a_1 \ge a_2 \ge \cdots \ge a_n$, so $a_1 < 1/2$ is the radius.

\subsection{Defender's Centroid Strategy} \label{sec:centroid-strat}
Consider a tree whose centroid breaks the tree into $n \geq 3$ components of lengths $a_1 \ge a_2 \ge \cdots \ge a_n$.
We introduce a strategy for the defender---called the {\em centroid} strategy---with which the defender hides the booby trap uniformly on the entire tree with probability 
\[
p = \frac{4 a_1}{1+4 a_1^2},
\]
and plants the booby trap at the centroid with probability 
\[
1- p = \frac{(1-2 a_1)^2}{1+4 a_1^2} \ge 0,
\]
so the centroid strategy is properly defined.
\begin{lemma} \label{lem:centroid-bound}
In the game $\Gamma'_1(Q)$ on a tree $Q$ whose centroid breaks the tree into $n \geq 3$ components of lengths $a_1 \ge a_2 \ge \cdots \ge a_n$, the defender's centroid strategy ensures that the expected payoff to the attacker is at most
\[
\frac{a_1}{1+4 a_1^2}.
\]
\end{lemma}
\textit{Proof.}
First, suppose the attacker chooses a subset of the tree that does not contain the centroid.
This subset must have length $x \le a_1 < 1/2$. The expected payoff is
\[
x(1- p  x) \le a_1 (1-p a_1)= \frac{a_1}{1+4 a_1^2},
\]
where the inequality follows from the fact that $x (1-p x)$ is increasing in $x$ for $x \le 1/2$.

Now suppose the attacker chooses a subset of the tree of length $x$ that contains the centroid. The expected payoff is
\[
p   (1-x) x \le p \cdot \frac{1}{4} = \frac{a_1}{1+4 a_1^2},
\]
and the proof is complete.
\hfill $\Box$

\bigskip

Note that as $a_1$ approaches $1/2$, the probability that the centroid strategy hides the booby trap at the centroid, namely $p$, approaches 0.
In the limit when $a_1 = 1/2$, we have $p=1$, and the centroid strategy reduces to the uniform strategy. 
In this case, the centroid strategy is optimal, by Lemma~\ref{lem:k=1}, because the tree can be partitioned into two components of length $1/2$.

\subsection{Attacker's Partition Strategy} \label{sec:partition}
Consider a tree whose centroid breaks the tree into $n \geq 3$ components of lengths $a_1 \ge a_2 \ge \cdots \ge a_n$.
One idea for the attacker is to treat the $n$ components as $n$ separate \textit{prize boxes} and ignore the centroid altogether.
In other words, the attacker can take the whole of any component---one prize box at a time---but never any subset that includes the centroid.
If the attacker chooses component $i$ with probability
\[
\frac{1/a_i}{\sum_{j=1}^n 1/a_j},
\]
for $i=1,\ldots,n$, then regardless of where the defender plants the booby trap, the attacker can guarantee an expected payoff of
\begin{equation}
\frac{n-1}{\sum_{j=1}^n 1/a_j}.
\label{eq:n-1}
\end{equation}

\cite{lidbetter2020search} considers the preceding discrete booby trap game with $n$ prize boxes---analogous to $n$ components in a tree after the centroid is removed---and 1 booby trap, and shows that the attacker can improve the expected payoff in \eqref{eq:n-1} by ignoring some boxes that contain the smallest prizes.
Intuitively, if $a_n$ is much smaller than $a_{n-1}$, then it is better for the attacker to ignore box $n$ entirely and just consider boxes $1, 2, \ldots, n-1$.
Let
\[
V(t) = \frac{t-1}{\sum_{j=1}^t 1/a_j},
\]
which is the expected payoff for the attacker, if the attacker chooses box $i$ with probability $(1/a_i) / \sum_{j=1}^t 1/a_j$, for $i=1,\ldots,t$ and ignores boxes $t+1, \ldots, n$.
\cite{lidbetter2020search} shows that the value of the discrete booby trap game is 
\begin{equation}
\max_{t=2,\ldots,n} V(t),
\label{eq:discrete_game}
\end{equation}
which is a lower bound for the expected payoff for the attacker on a tree whose centroid breaks the tree into $n \geq 3$ components having lengths $a_1 \ge a_2 \ge \cdots \ge a_n$.

Although the attacker in our booby trap game on a tree can guarantee the expected payoff in~\eqref{eq:discrete_game}, he can do better by using the centroid to connect some components.
For example, if instead of playing with components $1,\ldots, t$ to achieve $V(t)$, the attacker can use the centroid to connect components $t, t+1, \ldots, n$, so the last \textit{super component} now contains prize $\sum_{j=t}^n a_j$ rather than just $a_t$.
Write
\begin{equation}
U(t) = \frac{t-1}{\sum_{j=1}^{t-1} 1/a_j + 1/ \sum_{j=t}^n a_j} 
\label{eq:U(t)}
\end{equation}
for the expected payoff the attacker can guarantee in the discrete booby trap game with $t$ boxes containing $a_1, a_2, \ldots, a_{t-1}, \sum_{j=t}^n a_j$ with the equalizing strategy.
Because $\sum_{j=t}^n a_j > a_t$ for $t \le n-1$, it follows immediately that $U(t) > V(t)$, for $t=2, \ldots n-1$.


Generally speaking, for any subset $S \subseteq \{1,\ldots,n\}$, we can define an attacker's strategy as follows: take the whole of component $i \notin S$ with probability 
\[
\frac{1/ a_i}{\sum_{j \notin S} 1/a_j + 1/ \sum_{j \in S} a_j},
\]
or take all components in $S$ along with the centroid with probability
\[
\frac{1/ \sum_{j \in S} a_j}{\sum_{j \notin S} 1/a_j + 1/ \sum_{j \in S} a_j}.
\]
Letting $W(S)$ denote the expected payoff the attacker can guarantee with such a strategy, then
\begin{equation}
W(S) = \frac{n-|S|}{\sum_{j \notin S} 1/a_j + 1/ \sum_{j \in S} a_j}.
\label{eq:W(S)}
\end{equation}
We refer to the strategy in this class that maximizes the preceding as the attacker's \textit{partition strategy}.

What is the optimal choice of $S$ to maximize $W(S)$ in \eqref{eq:W(S)}?
Intuitively, it is better to consolidate smaller components so that the resulting \textit{prizes} in the \textit{boxes} are more balanced.
The next theorem formalizes this idea by showing that it is always best to consolidate components that are the smallest.

\begin{theorem}
\label{th:W(S)}
Consider the booby trap game $\Gamma'_1(Q)$ played on a tree $Q$ whose centroid breaks the tree into $n$ components of lengths $a_1 \ge a_2 \ge \cdots \ge a_n$.
To maximize $W(S)$ in \eqref{eq:W(S)}, it is sufficient to consider $S = \{t, t+1, \ldots, n\}$ for some $t = 2, 3, \ldots, n$.
In other words,
\[
\max_{S \subseteq \{1,2,\ldots,n\}} W(S) = \max_{t=2, \ldots, n} U(t),
\]
where $U(t)$ is given in \eqref{eq:U(t)}.
\end{theorem}
\textit{Proof.}
Consider two subsets of components, $S_1$ and $S_2$, which differ by exactly one component.
That is, there exists a subset $R$ such that $S_1 = R \cup \{j\}$ and $S_2 = R \cup \{k\}$, where $j \neq k$.
Without loss of generality, assume $j < k$ so $a_j \ge a_k$.
Our first step is to show that $W(S_1) \le W(S_2)$, which is equivalent to
\[
\frac{n - |S_1|}{\sum_{i \notin S_1} 1/a_i + 1/ \sum_{i \in S_1} a_i} \le \frac{n - |S_2|}{\sum_{i \notin S_2} 1/a_i + 1/ \sum_{i \in S_2} a_i}.
\]
Because $|S_1| = |S_2|$, the preceding---after writing $r = \sum_{i \in R} a_i$ for convenience---is equivalent to
\[
\frac{1}{a_k} + \frac{1}{r + a_j} \ge \frac{1}{a_j} + \frac{1}{r + a_k},
\]
or equivalently,
\[
\frac{a_j - a_k}{a_j a_k} \ge \frac{a_j-a_k}{(r+a_j) (r+a_k)}.
\]
This inequality is true because $a_j \ge a_k$ and $r > 0$, which shows that $W(S_1) \le W(S_2)$.

For any subset $S$ that is not in the form of $\{t, t+1, \ldots, n\}$ for some $t$, the attacker can improve (weakly) the expected payoff by replacing a component in $S$ with another component having a larger index.
Consequently, all subsets consisting of $n-t+1$ components are weakly dominated by $\{t, t+1, \ldots, n\}$, which completes the proof of this theorem.
\hfill $\Box$

\bigskip

From Theorem~\ref{th:W(S)}, to find the attacker's partition strategy, one needs to evaluate $U(t)$ in~\eqref{eq:U(t)} for $t=2,\ldots, n$.
It turns out that, we can cut the number of evaluations by a half, according to the next theorem.

\begin{theorem}
\label{th:U(t),n/2}
Consider the booby trap game $\Gamma'_1(Q)$ played on a tree $Q$ whose centroid breaks the tree into $n$ components of lengths $a_1 \ge a_2 \ge \cdots \ge a_n$.
To find the partition strategy for the attacker, it is sufficent to consider $U(t)$ for $t=2, \ldots, \lceil n/2 \rceil + 1$.
\end{theorem}
\textit{Proof.}
Our first task is to prove the claim that if $t > (n+1)/2$, then $U(t) > U(t+1)$, which is equivalent to
\[
\frac{t-1}{\sum_{j=1}^{t-1} 1/a_j + 1/ \sum_{j=t}^n a_j} > \frac{t}{\sum_{j=1}^{t} 1/a_j + 1/ \sum_{j=t+1}^n a_j}.
\]
Writing $r = \sum_{j=t+1}^n a_j$ for convenience, cross-multiplying and canceling common terms, the preceding is equivalent to
\[
\frac{t-1}{a_t} + \frac{t-1}{r} > \sum_{j=1}^{t-1} \frac{1}{a_j} + \frac{t}{a_t + r}.
\]
Because $a_1 \ge a_2 \ge \cdots \ge a_t$, it follows that
\[
\frac{t-1}{a_t} \ge \sum_{j=1}^{t-1} \frac{1}{a_j}. 
\]
If $t > (n+1)/2$, then
\[
\frac{t-1}{r} - \frac{t}{a_t+r} = \frac{1}{r (a_t + r)} ((t-1) a_t - r) \ge \frac{1}{r (a_t + r)} ((t-1) a_t - (n-t) a_t) > 0,
\]
where the first inequality uses the fact that $r = \sum_{j=t+1}^n a_j \le (n-t) a_t$ because $a_t \ge a_{t+1} \ge \cdots \ge a_n$.

In other words, if $n$ is even, then to maximize $U(t)$ we only need to consider $t$ up to $n/2+1$.
If $n$ is odd, then we need to consider $t$ up to $(n+1)/2 + 1$.
Combining the two cases gives the range specified in the theorem.
\hfill $\Box$

\bigskip

The results in Theorem~\ref{th:U(t),n/2} can be understood intuitively as follows.
For the case of $k=1$ booby trap, ideally the attacker likes to partition the network into 2 subnetworks each having length $1/2$, as shown in Lemma~\ref{lem:k=1}.
When using the centroid of a tree to connect components into the same subnetwork, intuitively it is better to include more components in the subnetwork as long as its length does not exceed $1/2$.
Because $a_1 \ge a_2 \ge \cdots \ge a_n$, it is clear that
\[
\sum_{j=1}^{\lceil n/2 \rceil} a_j \ge \frac{1}{2} \ge \sum_{j= \lceil n/2 \rceil +1}^n a_j.
\]
Therefore, it makes sense for the attacker to use the centroid to connect at least components $\lceil n/2 \rceil +1, \lceil n/2 \rceil +2, \ldots, n$ into a subnetwork, which is exactly what is stated in Theorem~\ref{th:U(t),n/2}.

\subsection{Bounds on the Value of the Game}
We now show that the ratio of the bounds given by the centroid strategy and the partition strategy is at most $27/25=1.08$. In other words, the defender's centroid strategy and the attacker's partition strategy both come within 8\% of the value of the game.
This implies we can bound the value of the game $\Gamma'_1(Q)$ on a tree $Q$ within 8\% of accuracy.

\begin{corollary} If $Q$ is a tree, the ratio of the expected payoff of the centroid strategy and the partition strategy in the game $\Gamma'_1(Q)$ is at most $27/25$.
\end{corollary} \label{cor:centroid-bound}
\textit{Proof.}
Suppose the removal of the centroid of $Q$ breaks the tree into $n \geq 3$ components of lengths $a_1 \ge a_2 \ge \cdots \ge a_n$. By Lemma~\ref{lem:centroid-bound} and Theorem~10, the ratio of the expected payoff of the centroid strategy and the partition strategy is
\[
\frac{a_1 / (1 +4 a_1^2)}{\max_{t=2,\ldots, \lceil n/2 \rceil +1} U(t)} \le \frac{a_1 / (1 +4 a_1^2)}{ U(2)} = \frac{a_1 / (1 +4 a_1^2)}{ a_1 (1-a_1)} = \frac{1}{(1-a_1) (1 + 4 a_1^2)}.
\]
Calculus shows that for $a_1 \in (0, 1/2)$, the preceding is maximized when $a_1=1/6$ and the maximized value is $27/25$.  \hfill $\Box$

\subsection{Centroid of a General Network}
We can extend the definition of a centroid to networks that are not trees. Unless a network $Q$ is 2-connected, it has at least one disconnecting point---a point whose removal disconnects the network. 
For such a point $x$, we extend the definition of $Q(x)$ and $h(x)$ in the natural way.
If $h(x) \le 1/2$, then we call $x$ a centroid of $Q$ and $h(x)$ its radius.
There cannot be more than one centroid, and the proof of this result is identical to the proof in Lemma~\ref{lem:centroid}.
In other words, each network either has a unique centroid or does not have one.
If a network has a centroid, then centroid strategy and the partition strategy can be defined identically to those in Subsections~\ref{sec:centroid-strat} and~\ref{sec:partition}, and the bounds given in Corollary~\ref{cor:centroid-bound} also hold.

\begin{figure}[!ht]
\center
\includegraphics[width=6cm]{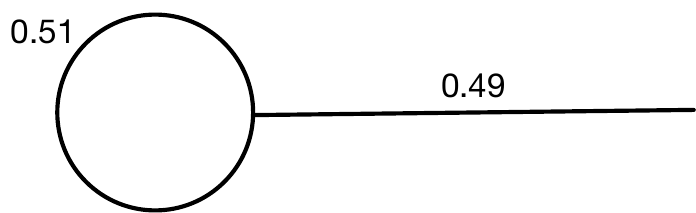}
\caption{A network that is not 2-connected, does not have a centroid, and can be partitioned into two connected subnetworks each having length $1/2$.}
\label{fig:ex1}
\end{figure}

Not all networks that are not 2-connected have a centroid.
One such example is a network that consists of a cycle of length $0.51$ and a single arc that has one endpoint connecting to the circle, as shown in Figure~\ref{fig:ex1}.
Removing any point on the circle (other than the node connecting the arc) does not disconnect the network.
Removing any point on the arc (including the node connecting the circle) disconnects the network into two components with one component having length at least $0.51 > 1/2$.
Therefore, this network is not 2-connected and it does not have a centroid.
However, this network can be partitioned into two connected subnetworks each having length $1/2$, so the value of the booby trap game with $k=1$ booby trap is $1/4$ according to Lemma~\ref{lem:k=1}.

\begin{figure}[!ht]
\center
\includegraphics[width=5cm]{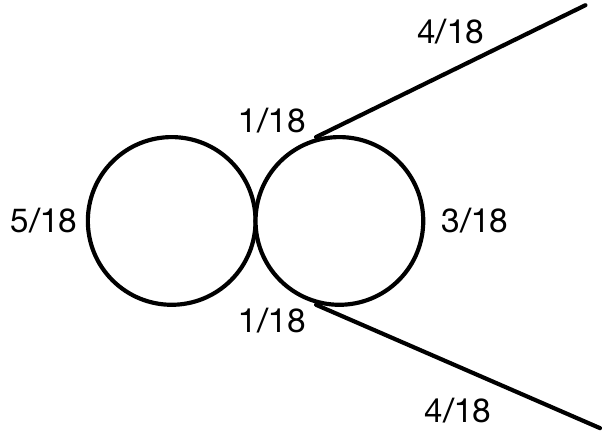}
\caption{A network that is not 2-connected, does not have a centroid, and cannot be partitioned into 2 subnetworks each having length $1/2$.}
\label{fig:ex2}
\end{figure}

Figure~\ref{fig:ex2} displays another network that is not 2-connected and does not have a centroid.
This network contains three nodes.
If we remove the node that connects an arc of length $4/18$ to the rest of the network, then the largest subnetwork has length $14/18 > 1/2$.
Removing the node connecting the two circles leaves the largest component having length $13/18 > 1/2$.
In addition, this network cannot be divided into two connected subnetworks each having length $1/2$, so Lemma~\ref{lem:k=1} does not apply.
A booby trap game played on this network requires further investigation with a careful examination of the network structure.

\section{Booby Trap Game on a Star Network}
\label{sec:star}
A star network consists of a node at the center where a number of arcs meet.
A star network is a tree, so we can use the results in Section~\ref{sec:tree} to compute an upper bound and a lower bound for the value of the game.
This section presents a few star networks for which we can determine the optimal strategy for each player and compute the value of the game.

\subsection{Symmetric Star}
Consider the booby trap game played on a symmetric star network with $k=1$ booby trap, which consists of $n$ arcs each having length $1/n$, for $n \ge 3$.
Theorem~\ref{thm:even-star} presents the solution if $n$ is even, and Theorem~\ref{thm:odd-star} presents the solution if $n$ is odd.

\begin{theorem} \label{thm:even-star}
Consider a symmetric star with $n$ arcs each having length $1/n$.
If $n$ is an even number, then the value of the game is $n/ (n^2+4)$.
The attacker's partition strategy is optimal, and the defender's centroid strategy is optimal.
\end{theorem}
\textit{Proof.}
First, evaluate the attacker's partition strategy.
From \eqref{eq:U(t)}, the attacker's partition strategy guarantees expected payoff
\[
\max_{t=2,\ldots, \lceil n/2 \rceil + 1} U(t) = \frac{t-1}{(t-1) n + n/ (n-t+1)}. 
\]
Taking $t = n/2 + 1$ yields
\[
U \left( \frac{n}{2} + 1 \right) = \frac{n}{n^2+4},
\]
which is a lower bound for the value of the game.
In other words, the attacker can guarantee the expected payoff $n/ (n^2+4)$ by taking either a group of $n/2$ arcs with probability $4/(n^2 + 4)$, or each of the remaining arcs with probability $2n/ (n^2 + 4)$.


Second, evaluate the defender's strategy.
By Lemma~\ref{lem:centroid-bound}, the defender's centroid strategy guarantee's the attacker's payoff at most
\[
\frac{a_1}{1 +4 a_1^2} = \frac{1/n}{1 +4/n^2} = \frac{n}{n^4+4},
\]
which is an upper bound for the value of the game.
Because the attacker's partition strategy with $t=n/2+1$ and the defender's centroid strategy both achieve the same value, each strategy is optimal and the achieved value is the value of the game.
\hfill $\Box$

\bigskip

The results in Theorem~\ref{thm:even-star} can be extended immediately to networks that have a centroid, which breaks the network into an even number of subnetworks each having the same length, which is stated below.

\begin{corollary} \label{co:even-symmetric}
Consider a network that has a centroid, and the centroid breaks the network into $n$ subnetworks each having length $1/n$.
If $n$ is an even number, then the value of the game is $n/ (n^2+4)$.
The attacker's partition strategy is optimal, and the defender's centroid strategy is optimal.
\end{corollary}

If a network can be partitioned into 2 subnetworks each having length $1/2$, then the attacker's strategy that takes each subnetwork with probability $1/2$ is optimal.
Based on that observation, it is somewhat intuitive that in a symmetric star, a partition strategy that consolidate half of the arcs into a single subnetwork would work well.
Theorem~\ref{thm:even-star} shows that such partition strategy is indeed optimal for the attacker for a symmetric star with an even number of arcs.

When the number of arcs is odd, then it turns out that it is optimal for the attacker to use a partition strategy by consolidating either $(n-1)/2$ or $(n+1)/2$ arcs into the same subnetwork.
The defender's centroid strategy, however, is no longer optimal.
When not planting the booby trap at the center node, instead of planting the booby trap uniformly randomly, the defender can guarantee a smaller expected payoff for the attacker by planting the booby trap at places nearer the center node but not near the end of an arc.

%

\bigskip



\begin{figure}[ht]
\center
\includegraphics[width=6cm]{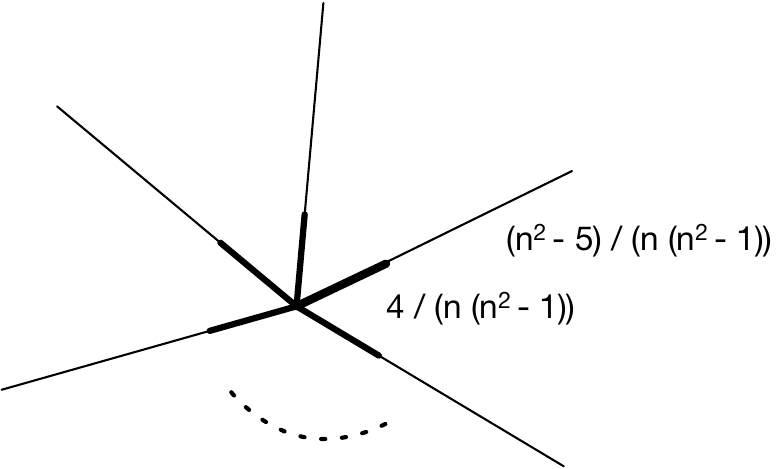}
\caption{The defender's optimal strategy for a symmetric star network with $n$ arcs when $n$ is odd: Plant the booby trap in the center with probability $(n^2-4n+3)/(n^2+3)$, or plant it uniformly randomly in the subset within distance $4/ (n (n^2-1))$ from the center.}
\label{fig:symmetric_odd}
\end{figure}

\begin{theorem} \label{thm:odd-star}
Consider a symmetric star with $n$ arcs each having length $1/n$.
If $n$ is an odd number, then the value of the game is $(n^2-1)/ (n(n^2+3))$.
The attacker's partition strategy is optimal.
An optimal strategy for the defender is to place the booby trap at the center with probability $(n^2-4n+3)/(n^2+3)$; or place it uniformly randomly in the subset consisting of all points within distance $4/(n (n^2-1))$ from the center with probability $4n / (n^2+3)$, as shown in Figure~\ref{fig:symmetric_odd}.

\end{theorem}
\textit{Proof.}
First, evaluate the attacker's partition strategy in
\eqref{eq:U(t)} by taking $t = (n+3)/2$.
In other words, the attacker either takes a group of $(n-1)/2$ arcs or any one of the remaining $(n+1)/2$ arcs.
With this strategy, the attacker guarantees expected payoff
\[
U\left(\frac{n+3}{2} \right) = \frac{(n+1)/2}{n ( n+1)/2 + 1 / ((n-1)/(2n))} = \frac{n^2-1}{n (n^2+3)}.
\]
One can verify that taking $t= (n+1)/2$ also achieves this same expected payoff, which is a lower bound for the value of the game.

\bigskip

Second, consider the defender's strategy.
Suppose first that the attacker chooses a subset of length $x \leq 1/n$.
It is clear that the attacker should avoid the center and choose a subinterval of an arc starting from the far end to minimize the probability of running into the booby trap.
By choosing $x = (n^2-5) / (n (n^2-1))$, the attacker guarantees success because the defender never places the booby trap in the subinterval of length $(n^2-5)/(n (n^2-1))$ at the far end of each arc, as seen in Figure~\ref{fig:symmetric_odd}.
For $x \ge (n^2-5) / (n (n^2-1))$, the probability of running into the booby trap is
\[
\frac{x -  (n^2-5) / (n (n^2-1))}{4 n / (n (n^2-1))} \frac{4n}{n^2+3} = \left( x - \frac{n^2-5}{n (n^2-1)} \right) \frac{n (n^2-1)}{n^2+3},
\]
so the expected payoff is
\[
x \cdot \left( 1 - \left( x - \frac{n^2-5}{n (n^2-1)} \right) \frac{n (n^2-1)}{n^2+3} \right) = - \frac{n (n^2-1)}{n^2+3} x^2 + \frac{2(n^2-1)}{n^2+3} x.
\]
Calculus shows that the preceding quadratic function increases in $x$ for $x \in [(n^2-5) / (n (n^2-1)), 1/n]$.
When $x=1/n$, the maximal expected payoff is $(n^2-1)/ (n (n^2+3))$.

Now suppose that the attacker chooses a subset of length $x > 1/n$.
The subset must include the center node.
In addition, the best subset to minimize the probability of containing the booby trap must consist of $j = \lfloor x/ (1/n) \rfloor = \lfloor n x \rfloor \ge 1$ whole arcs and an interval of length $y = x- \lfloor n x \rfloor < 1/n$ of another arc spanning from the center.
If $y > 4 / (n (n^2-1))$ then it is best to take $y \rightarrow 1/n$, or equivalently, $j+1$ whole arcs.
Hence, it is sufficient to consider the case $y \leq 4 / (n (n^2-1))$, whose expected expected payoff is
\begin{align*}
&\left( \frac{j}{n} + y \right) \left( (n-j-1) \frac{4}{(n^2+3)} + \frac{4}{(n^2+3)} \left( 1 - \frac{y}{4/(n(n^2-1))} \right) \right) \\
&= \left( \frac{j}{n} + y \right) \left(  \frac{4 (n-j)}{(n^2+3)} - \frac{n (n^2-1) y}{(n^2+3)} \right) \\
&= - \frac{n (n^2-1)}{n^2+3} y^2 + \left( \frac{4n}{n^2+3} -j \right) y + \frac{4j(n-j)}{n (n^2+3)}.
\end{align*}
Because $4n / (n^2+3) -j \le 0$ for $n \ge 3$ and $j \ge 1$, the preceding quadratic function decreases in $y$ for $y \geq 0$.
Therefore, it is best to take $y=0$.
In other words, it is best to take $j$ whole arcs to obtain expected payoff
\[
\frac{4j(n-j)}{n (n^2+3)}
\]
for $j=1, \ldots, n$.
To maximize the preceding, it is optimal to take $j= (n-1)/2$ or $(n+1)/2$, each of which results in an expected payoff
\[
\frac{n^2-1}{n (n^2+3)}.
\]
Consequently, we have shown that the defender's strategy guarantees the attacker's payoff at most $(n^2-1)/(n(n^2+3))$, which is therefore an upper bound for the value of the game.

The theorem follows because the lower bound for the value of the game obtained by the attacker's strategy coincide with the upper bound for the value of the game obtained by the defender's strategy.
\hfill $\Box$



\bigskip

Using the results in Theorems~\ref{thm:even-star} and \ref{thm:odd-star}, it is straightforward to verify that the value of the game for the booby trap game played on a symmetric star network decreases in $n$---the number of arcs---and asymptotically approaches $1/n$ as $n$ tends to infinity.

\subsection{Star with Three Arcs}
This section presents the case, in which the space $Q$ is a star network with three arcs and the defender has $k=1$ booby trap.
Write $a_j$ for the length of arc $j$, for $j=1,2,3$, with $a_1 \ge a_2 \ge a_3$, without loss of generality.
If $a_1 \ge 1/2$, then it follows from Lemma~\ref{lem:k=1} that the value of the game is $1/4$, because the network can be partitioned into two connected subsets of equal measure.
The next theorem solves the game when $a_1 < 1/2$.

\begin{theorem}
Consider the booby trap game played on a star with 3 arcs and 1 booby trap, with arc lengths $a_1 \ge a_2 \ge a_3$.
If $a_1 < 1/2$, the value is $a_1(1-a_1)$. 
The attacker's partition strategy is optimal.
For the defender, it is optimal to place the booby trap on arc $j$ at distance $x$ from the centroid, where $x$ is chosen according to the cumulative distribution function $F_j(x) = \min \{a_j, (1-a_j) x/a_j \}$.
\end{theorem}
\textit{Proof.}
Consider the attacker's partition strategy with $t=2$, which guarantees an expected payoff
\[
U(2)= \frac{2-1}{1/a_1 + 1/(a_2 + a_3)} = a_1 (1-a_1),
\]
because $a_2 + a_3 = 1 - a_1$.
Hence, the preceding is a lower bound for the value of the game.


Consider the defender's strategy.
An equivalent description of the attacker's strategy is to place the booby trap on arc $i$ with probability $a_i$, for $i=1,2,3$.
If arc $i$ is chosen, for $i=1,2,3$, then the booby trap is placed uniformly randomly on the segment starting at the centroid with length $a_i^2 / (1-a_i)$, which is strictly less than $a_i$ because $a_i < 1/2$.
Consequently, for each arc $i$ there is zero probability that the booby trap is at a distance greater than $a_i^2/(1-a_i)$ on that arc, $i=1,2,3$.

First, suppose that the attacker chooses a subset of arc $j$, for some $j=1,2,3$. Clearly, the attacker should choose a subset consisting of all points on that arc, whose distance from the centroid is greater than some $x \le a_i^2/(1-a_i)$. In this case, the expected payoff is
\[
\left( \frac{(1-a_i)x}{a_i} + 1-a_i \right)(a_i-x) = a_i(1-a_i) - \frac{(1 - a_i)x^2}{a_i} \le a_1(1-a_1),
\]
because $a_1$ maximizes $a_i (1-a_i)$.

Second, suppose that the attacker chooses a subset that has nonempty intersection with exactly two arcs. The expected payoff is maximized if the attacker chooses the whole of one arc $i$ along with all the points on some other arc $j \ne i$ within some distance $x \le a_j^2/(1-a_j)$ of the centroid, which is equal to
\begin{align}
(a_i+x) \left( 1-a_i- \frac{(1-a_j)x}{a_j} \right).
\label{eq:n=3}
\end{align}
Because there are only 3 arcs, either $i=1$ or $j=1$.
Consider the two cases separately:
\begin{enumerate}
\item 
If $i=1$, then the expected payoff in \eqref{eq:n=3} reduces to
\[
 a_1(1-a_1) -\frac{(1-a_j)x^2}{a_j} - \frac{(a_1-a_j)x}{a_j}  \le a_1(1-a_1),
\]
because $a_1 \ge a_j$.
\item
If $j=1$, then the expected payoff in \eqref{eq:n=3} becomes $(a_i+x)(1-a_i-(1-a_1)x/a_1)$, which is a quadratic function in $x$.
The quadratic function is maximized at $x=(1-2a_i+a_1a_i)/(2(1-a_1))$, where the expected payoff is
\begin{align*}
\frac{(a_1+a_i-2a_1a_i)^2}{4a_1(1-a_1)} = a_1(1-a_1)\left( \frac{a_i(1-2a_1)+a_1}{a_1(1-2a_1)+a_1} \right)^2  \le a_1(1-a_1),
\end{align*}
since $a_i \le a_1$ and $1-2a_1 \ge 1$.
\end{enumerate}

Finally, suppose the attacker chooses a subset that has nonempty intersection with all three arcs. The expected payoff is maximized if the attacker chooses all the points on some arc $i$ with some distance $x \le a_i^2/(1-a_i)$ of the centroid along with the whole of the other two arcs. The expected payoff is
\begin{align*}
(1-(a_i-x))\left(a_i - \frac{(1-a_i)x}{a_i}\right)  = a_i(1-a_i) - \frac{x^2}{a_i}-\frac{(1-2a_i)x}{a_i}  \le a_1(1-a_1),
\end{align*}
because $a_i \le 1/2$ and $a_1$ maximizes $a_i(1-a_i)$.

In summary, we have shown that the defender's strategy guarantees the attacker's expected payoff at most $a_1 (1-a_1)$, which is therefore an upper bound for the value of the game.
The theorem then follows because the lower bound for the value of the game obtained by the
attacker’s strategy coincide with the upper bound for the value of the game obtained by the defender’s strategy
\hfill $\Box$

\section{Conclusions}
\label{sec:conclusions}
We have introduced the study of continuous booby trap games, played either on a subset of Euclidean space or on a network. We gave a full solution of the game for the first model, and for the second model we solved the game in a number of cases. In the case of one booby trap on a tree network, we defined strategies for both players that are close to optimal. In fact, the attacker's partition strategy is optimal in all the examples of star network for which we have a complete solution, and we conjecture that this strategy is optimal for all star networks, and even for all trees.

To the best of our knowledge, this paper is the first to address a booby trap game played in a continuous space.
Although this paper presents many interesting results, solving the booby trap game in general appears to be quite challenging.
Immediate future research directions include developing strong heuristic strategies useful in practice and algorithms that approximate optimal strategies for general cases.

\bibliographystyle{apalike}
\bibliography{ref}

\end{document}